\documentclass[12pt]{article}

\usepackage{amsmath,amsthm}

\usepackage{amsfonts}
\usepackage{natbib}

\usepackage[dvips]{graphicx}

\makeatletter
\newcommand{\@makemycaption}[2]{%
\vspace{10pt}%
{\textbf{#1}:#2\par}%
}
\renewcommand{\figure}{%
\let\@makecaption\@makemycaption\@float{figure}}
\renewcommand{\table}{%
\let\@makecaption\@makemycaption\@float{table}}
\makeatother

\date{}

\newcommand{\eps}{\varepsilon}
\newcommand{\To}{\rightarrow}
\newcommand{\altild}{\widetilde{\alpha}}

\numberwithin{equation}{section}

\begin{document}

\begin{center}
  \textbf{\large Generalized Likelihood Ratio Statistics and
    Uncertainty Adjustments in Efficient Adaptive Design of Clinical
    Trials}

\bigskip

\textbf{Jay Bartroff}\\ Department of Mathematics, University of Southern California, Los
 Angeles, CA, USA {\renewcommand{\thefootnote}{}\footnote{Address
 correspondence to Jay Bartroff, Department of Mathematics, University of Southern California, 3620 South Vermont Ave., KAP 108, Los Angeles, CA 90089-2532, USA; E-mail: bartroff@usc.edu}} 

\medskip

\textbf{Tze Leung Lai}

Department of Statistics, Stanford University, Stanford, CA, USA
\end{center}

\bigskip

\noindent\textbf{\small Abstract:} {\small A new approach to adaptive design of clinical trials is proposed in a general
  multiparameter exponential family setting, based on generalized
  likelihood ratio statistics and optimal sequential testing theory.
  These designs are easy to implement, maintain the prescribed Type~I
  error probability, and are asymptotically efficient.  Practical
  issues involved in clinical trials allowing mid-course adaptation and the large literature
  on this subject are discussed, and comparisons between the proposed
  and existing designs are presented in extensive simulation studies
  of their finite-sample performance, measured in terms of the
  expected sample size and power functions.}

\vspace{1.5in}

\noindent{\small\textbf{Keywords.} Hoeffding's
information bound; Internal pilot; Kullback-Leibler information; Modified
Haybittle-Peto test; Multiparameter exponential family; Sample size re-estimation}

\noindent{\small\textbf{Subject Classifications.} 62L10; 62F03; 62P10}

\newpage

\section{\large INTRODUCTION}\label{sec:intro}

\noindent Because of the ethical and economic considerations in the
design of clinical trials to test the efficacy of new treatments and
because of lack of information on the magnitude and sampling
variability of the treatment effect at the design stage, there has
been increasing interest from the biopharmaceutical industry in
sequential methods that can adapt to information acquired during the
course of the trial. Beginning with~\citet{Bauer89}, who introduced
sequential adaptive test strategies over a planned series of separate
trials, and~\citet{Wittes90}, who discussed internal pilot studies, a large literature has grown on adaptive
design of clinical trials. Depending on the topics covered, the term
``adaptive design'' in this literature is sometimes replaced by
``sample size re-estimation,'' ``trial extension'' or ``internal
pilot studies.''  In standard clinical trial designs, the sample size
is determined by the power at a given alternative, but in practice, it
is often difficult for investigators to specify a realistic
alternative at which sample size determination can be based. Although
a standard method to address this difficulty is to carry out a
preliminary pilot study, the results from a small pilot study may be
difficult to interpret and apply, as pointed out by~\citet{Wittes90},
who proposed to treat the first stage of a two-stage clinical trial as
an internal pilot from which the overall sample size can be
re-estimated. The problem of sample size re-estimation based on
observed treatment difference at some time before the prescheduled end
of a clinical trial has attracted considerable attention during the
past decade; see, e.g., \citet{Gould92}, \citet{Herson93},
\citet{Birkett94}, \citet{Denne99,Denne00},
\citet[Section~14.2]{Jennison00}, \citet{Shih01} and
\citet{Whitehead01}. For normally distributed outcome variables with
known variances, \citet{Proschan95}, \citet{Fisher98}, \citet{Posch99}
and~\citet{Shen99} have proposed ways to adjust the test statistics
after mid-course sample size modification so that the Type I error
probability is maintained at the prescribed level. By making use of a
generalization of the Neyman-Pearson lemma, \citet{Tsiatis03} showed
that these adaptive tests of a simple null versus a simple alternative
hypothesis are inefficient because they are not based on likelihood
ratio statistics. \citet{Jennison03} gave a general weighted form of
these test statistics and demonstrated in simulation studies that the
adaptive tests performed considerably worse than group sequential
tests. \citet{Jennison06a} recently introduced adaptive group
sequential tests that choose the $j$th group size and stopping
boundary on the basis of the cumulative sample size $n_{j-1}$ and the
sample sum $S_{n_{j-1}}$ over the first $j-1$ groups, and that are
optimal in the sense of minimizing a weighted average of the expected
sample sizes over a collection of parameter values subject to
prescribed error probabilities at the null and a given alternative
hypothesis. They showed how the corresponding optimization
problem can be solved numerically by using the backward induction
algorithms for ``optimal sequentially planned'' designs developed by
\citet{Schmitz93}. \citet{Jennison06b} found that standard (non-adaptive) group
sequential tests with the first stage chosen optimally are nearly as
efficient as their optimal adaptive tests.

Except for Jennison and Turnbull's optimal adaptive group sequential
tests and the extensions of the sample size re-estimation approach to
group sequential testing considered by~\citet{Cui99},
\citet{Lehmacher99} and~\citet{Denne00}, previous works in the
literature on mid-course sample size re-estimation have focused on
two-stage designs whose second-stage sample size is determined by the
results from the first stage (internal pilot), following the seminal
work of~\citet{Stein45} in this area. Although this approach is
intuitively appealing, it does not adjust for the uncertainty in the
first-stage parameter estimates that are used to determine the
second-stage sample size. Moreover, it considers primarily the special
problem of comparing the means of the two normal populations, using
the central limit theorem for extensions to more general situations.
The case of unknown common variance at a prespecified alternative for
the mean difference was considered first, as in~\citet{Stein45} and
the first set of references in the preceding paragraph. Then the case
of known variances in the absence of a prespecified alternative for
the mean difference was studied, as in the second set of references
above. 

\citet{Bartroff08c} recently gave a unified treatment of both
cases in the general framework of multiparameter exponential families.
It uses efficient generalized likelihood ratio (GLR) statistics in
this framework and adds a third stage to adjust for the sampling
variability of the first-stage parameter estimates that determine the
second-stage sample size. Specifically, let
$X_1,X_2,\ldots$ be independent $d$-dimensional random vectors from a
multiparameter exponential family $f_\theta (x) = \exp\{\theta' x -
\psi(\theta)\}$ of densities with respect to some measure $\nu$ on
${\bf R}^d$. Let $S_n = X_1 + \ldots + X_n$. A sufficient statistic
based on $(X_1, \ldots, X_n)$ is the sample mean $\bar X_n = S_n / n$,
which is the maximum likelihood estimate of the mean
$\nabla\psi(\theta)$. Consider the hypothesis $u (\theta) = u_j$,
where $j = 0$ or $1$, $u : \Theta \rightarrow {\bf R}$ is continuously
differentiable and $\Theta = \{\theta : \int e^{\theta' x} d \nu (x)
< \infty\}$ is the natural parameter space. As noted by
\citet[p.~513]{Lai04}, the GLR statistic
$\Lambda_{i, j}$ for total sample size $n_i$ at stage $i$ has the form
\begin{equation}\label{eq:1}\Lambda_{i, j} = n_i \{\widehat{\theta}_{n_i}' \bar X_{n_i} - \psi(\widehat
\theta_{n_i})\} - \sup_{\theta : u(\theta) = u_j} n_i \{\theta' \bar
X_{n_i} - \psi(\theta)\} = \inf_{\theta : u(\theta) = u_j} n_i
I(\widehat \theta_{n_i}, \theta) ,\end{equation} in which
$\widehat{\theta}_{n}=(\nabla\psi)^{-1}(\overline{X}_{n})$ and $I(\theta, \lambda)$ is the Kullback-Leibler information
number given by \begin{equation}\label{eq:2}I(\theta, \lambda) = E_\theta [\log \{f_\theta (X_i) / f_\lambda (X_i)\}]
= (\theta - \lambda)' \nabla\psi(\theta) - \{\psi(\theta) - \psi
(\lambda)\}. \end{equation} The
possibility of adding a third stage to improve two-stage designs dated
back to~\citet{Lorden83}. Whereas Lorden used crude upper bounds for the Type I error
probability that are too conservative for practical applications,
\citet{Bartroff08c} overcome this difficulty by developing numerical methods to compute
the Type I error probability, and also extended
the three-stage test to multiparameter and
multi-armed settings, thus greatly broadening the scope
of these efficient adaptive designs. A review of their method is given
in Section~\ref{sec:three}. In Section~\ref{sec:asymptotics} we prove
the asymptotic optimality of these adaptive tests in the
multiparameter case, extending Lorden's~(\citeyear{Lorden83}) result
for the special case $d=1$. 

Another new addition to this asymptotic optimality theory of adaptive
designs is related to the problem of trial extension considered in
Section~\ref{sec:Mtild}. As pointed out by~\citet{Cui99}, the issue of increasing
the maximum sample size sometimes arises after interim analysis in
group sequential trials. They cited a study protocol, which was
reviewed by the Food and Drug Administration, involving a Phase~III
group sequential trial for evaluating the efficacy of a new drug to
prevent myocardial infarction in patients undergoing coronary artery
bypass graft surgery.  During interim analysis, the observed incidence
for the drug achieved a reduction that was only half of the target
reduction assumed in the calculation of the maximum sample size $M$,
resulting in a proposal to increase the maximum sample size to
$\widetilde{M}$ ($N_{\max}$ in their notation). \citet{Cui99} and
\citet{Lehmacher99} extended the sample size re-estimation approach to
adaptive group sequential trials by adjusting the test statistics as
in~\citet{Proschan95} and allowing the future group sizes to be
increased or decreased during  interim analyses so that the overall
sample size does not exceed $\widetilde{M}(>M)$ and the Type I error
probability is maintained at the prescribed level. In
Section~\ref{sec:Mtild} we propose an alternative approach that is
shown to be asymptotically efficient in
Section~\ref{sec:asymptotics}. Whereas the adaptive designs in
Section~\ref{sec:three} assume a given maximum sample size $M$ and
require at most 3 stages, Section~\ref{sec:Mtild} extends them to allow
mid-course increase of $M$ to $\widetilde{M}$. These adaptive designs
involve no more than 4 stages and an adjustment in the maximum sample
size is made at the third stage. Computational algorithms to implement
them are provided in Section~\ref{sec:implementation} and simulation
results on their performance are given in Section~\ref{sec:Mtild_eg}.

\citet{Bartroff08c} carried out comprehensive simulation studies of the
performance, measured in terms of the expected sample size and power
functions, of their adaptive test and compared it with other adaptive
tests in the literature.  In the case of normal mean with known
variance and Type~I and II error constraints under the null and a
given alternative hypothesis, they showed that their adaptive test is
comparable to the benchmark optimal adaptive test of
\citet{Jennison06a,Jennison06b}, which is superior to the existing
two-stage adaptive designs. On the other hand, whereas the benchmark
optimal adaptive test needs to assume a specified alternative, these
adaptive two-stage tests and the adaptive tests of~\citet{Bartroff08c}
do not require such assumptions as they consider the estimated
alternative at the end of the first stage. In their recent survey of
adaptive designs, \citet{Burman06} pointed out that previous
criticisms of the statistical properties of two-stage adaptive designs
may be unconvincing in some situations when flexibility and not having
to specify parameters that are unknown at the beginning of a trial
(like the relevant treatment effect or variance) are more imperative
than efficiency or being powerful. The adaptive designs in
\citet{Bartroff08c} and this paper can fulfill the seemingly disparate
requirements of flexibility and efficiency on a design. Rather than
achieving exact optimality at a specified collection of alternatives
through dynamic programming, they achieve asymptotic optimality over
the entire range of alternatives, resulting in near-optimality in
practice. They are based on efficient test statistics of the GLR type,
which have an intuitively ``adaptive'' appeal via estimation of
unknown parameters by maximum likelihood, ease of implementation and
freedom from having to specify the relevant alternative; see
Section~\ref{sec:design}.

\citet{Bauer06} have found from a search of the medical literature
that adaptive designs have not been widely used in practice and that
``adaptations in practice are rather limited to sample size
reassessment.'' Perhaps one reason why these two-stage adaptive
designs have not gained wide acceptance is their use of seemingly
unnatural and convoluted test statistics (e.g., the inefficient test
statistics mentioned in the first paragraph). This can be circumvented
by the use of efficient GLR statistics in our adaptive tests with no
more than 3 (or 4) stages. Another reason may be the lack of routine
medical studies to which adaptive designs can lead to substantial
improvements over current practice. In Section~\ref{sec:phaseII} we
consider one such potential application in Phase~II cancer studies. We
show how our adaptive designs offer improvements over
Simon's~(\citeyear{Simon89}) optimal two-stage designs, which are
commonly used in single-arm cancer trials, and over their analogs, due
to~\citet{Thall88}, for randomized trials in
Section~\ref{sec:phaseII-III}. Section~\ref{sec:var_unknown} considers
another issue that often arises in the design of clinical trials,
namely, nuisance parameters. ``Very often, statistical information
also depends on nuisance parameters (e.g., the standard deviation of
the response variable). Extension of statistical information is a
design adaptation that occurs most frequently''~\citep{Hung06}. The
simulation results in Section~\ref{sec:var_unknown} shows how our
adaptive test resolves the difficulties with conventional two-stage
designs to treat nuisance parameters, which are estimated at the end
of the first stage and then used to estimate the second stage sample
size.  Further discussion of our proposed approach to adaptive designs
and some concluding remarks are given in Section~\ref{sec:discussion}.

\section{\large EFFICIENT ADAPTIVE DESIGN AND GLR TESTS}\label{sec:design}
\subsection{An Adaptive 3-Stage GLR Test}\label{sec:three}

Whereas~\citet{Tsiatis03} consider the case of simple null and
alternative hypotheses $\theta=\theta_j$ ($j=0,1$) for which
likelihood ratio tests are most powerful even in their group
sequential designs, \citet{Bartroff08c} use the GLR
statistics~(\ref{eq:1}) in an adaptive three-stage test of the
composite null hypothesis $H_0: u(\theta)\le u_0$, where $u$ is a
smooth real-valued function such that
\begin{equation}\label{eq:30}\mbox{$I(\theta,\lambda)$ is increasing
    in $u(\lambda)$ for every fixed $\theta$.}\end{equation} Let $n_1=m$ be the sample size of the
first stage (or internal pilot study) and $n_3=M$ be the maximum total
sample size, both specified before the trial. Let $u_1>u_0$ be the
alternative implied by the maximum sample size $M$ and the reference
Type II error probability $\altild$.  That is, $u_1(>u_0)$ is the
alternative where the fixed sample size (FSS) GLR test with Type I
error probability $\alpha$ and sample size $M$ has power
$\inf_{\theta:u(\theta)=u_1} P_\theta\{\mbox{Reject $H_0$}\}$ equal to $1-\altild$, as in~\citet[Section~3.4]{Lai04}. The
three-stage test of $H_0: u(\theta)\le u_0$ stops and rejects $H_0$ at
stage $i\le 2$ if \begin{equation}\label{eq:6}n_i<M,\quad u(\widehat
  \theta_{n_i}) > u_0\quad\mbox{and}\quad\Lambda_{i, 0} \geq
  b.\end{equation} Early stopping for futility (accepting $H_0$) can
also occur at stage $i \leq 2$ if
\begin{equation}\label{eq:7}n_i<M,\quad u(\widehat \theta_{n_i}) <
u_1\quad\mbox{and}\quad\Lambda_{i, 1} \geq \widetilde
b. \end{equation} The test rejects $H_0$ at stage $i=2$ or $3$ if
\begin{equation}\label{eq:8}n_i=M,\quad
u(\widehat{\theta}_M)>u_0\quad\mbox{and}\quad \Lambda_{i,0}\ge
c,\end{equation} accepting $H_0$ otherwise. The sample size $n_2$ of
the three-stage test is given by \begin{equation}\label{eq:3}n_2
=m\vee \{M \wedge \lceil (1 + \rho_m) n(\widehat
\theta_m)\rceil\}\end{equation} with
\begin{equation}\label{eq:5}n(\theta) = \min\{|\log \alpha | /
\inf_{\lambda : u (\lambda) = u_0} I(\theta, \lambda), |\log
\widetilde \alpha | / \inf_{\lambda : u(\lambda) = u_1} I(\theta,
\lambda)\},\end{equation} where $\rho_m>0$ is an inflation factor to
adjust for uncertainty in $\widetilde{\theta}_m$; see the examples in
Section~\ref{sec:examples}. Letting $0<\eps,\widetilde{\eps}<1$,
define the thresholds $b, \widetilde{b}$ and $c$ to satisfy the
equations
\begin{eqnarray}&&\sup_{\theta:u(\theta)=u_1} P_\theta
  \{\mbox{(\ref{eq:7}) occurs for $i = 1$ or $2$}\} =
\widetilde{\eps} \altild,\label{eq:9}\\
&&\sup_{\theta:u(\theta)=u_0} P_\theta \{\mbox{(\ref{eq:7}) does not
  occur for $i\le 2$, (\ref{eq:6}) occurs for $i=1$ or $2$}\} = \eps \alpha ,\label{eq:10}\\
&& \sup_{\theta:u(\theta)=u_0} P_\theta \{\mbox{(\ref{eq:6}) and
  (\ref{eq:7}) do not occur for $i\le 2$, (\ref{eq:8}) occurs}\} =
(1-\eps)\alpha. \label{eq:11}\end{eqnarray} The probabilities in
(\ref{eq:9})-(\ref{eq:11}) can be computed by using the normal approximation to the signed-root
likelihood ratio statistic
$$\ell_{i,j}=\{\mbox{sign}(u(\widehat{\theta}_{n_i})-u_j)\}
(2n_i\Lambda_{i,j})^{1/2},$$
($1\le i\le 3$ and $j=0,1$) under
$u(\theta)=u_j$, as in~\citet[p.~513]{Lai04}. When
$u(\theta)=u_j$, $\ell_{i,j}$ is approximately normal with mean 0, variance
$n_i$, and the increments $\ell_{i,j}-\ell_{i-1,j}$ are asymptotically
independent.  We can therefore approximate $\ell_{i,j}$ by a sum of independent standard normal random variables under
$u(\theta)=u_j$ and thereby determine $b, \widetilde{b}$ and $c$.
Note that this normal approximation can also be used for the choice of
$u_1$ implied by $M$ and $\altild$. Computational details are given in
Section~\ref{sec:implementation}, as well as an alternate method for
computing the thresholds by Monte Carlo.

A special multiparameter case of particular interest in clinical
trials involves $K$ independent populations having density functions
$\exp\{\theta_k x-\widetilde{\psi}_k(\theta_k)\}$ so that $\theta'
x-\psi(\theta)=\sum_{k=1}^K \{\theta_k
x_k-\widetilde{\psi}(\theta_k)\}$. In multi-armed trials, for which
different numbers of patients are assigned to different treatments,
the GLR statistic $\Lambda_{i,j}$ for testing the hypothesis
$u(\theta_1,\ldots,\theta_K)=u_j$ ($j=0$ or $1$) at stage $i$ has the
form $$\Lambda_{i,j} = \sum_{k=1}^K
n_{ki}\{\widehat{\theta}_{k,n_{ki}} \overline{X}_{k,n_{ki}}-
\widetilde{\psi}(\widehat{\theta}_{k,n_{ki}}) \}-\sup_{\theta:
  u(\theta_1,\ldots,\theta_K)=u_j} \sum_{k=1}^K n_{ki}\{\theta_k
\overline{X}_{k,n_{ki}}-\widetilde{\psi}(\theta_k)\},$$
in which
$n_{ki}$ is the total number of observations from the $k$th population
up to stage $i$. Let $n_i=\sum_{k=1}^K n_{ki}$. As pointed out in
Section~3.4 of~\citet{Lai04}, the normal approximation to the signed
root likelihood ratio statistic is still applicable when
$n_{ki}=p_kn_i+O_p(n_i^{\frac{1}{2}})$, where $p_1,\ldots,p_K$ are
nonnegative constants that sum up to 1, as in random allocation of
patients to the $K$ treatments (for which $p_k=1/K$).

\subsection{Mid-course Modification of Maximum Sample Size}\label{sec:Mtild}

We now modify the adaptive designs in the preceding section to
accommodate the possibility of mid-course increase of the maximum
sample size from $M$ to $\widetilde{M}$. Let $u_2$
be the alternative implied by $\widetilde{M}$ so that the
level-$\alpha$ GLR test with sample size $\widetilde{M}$
has power $1-\widetilde{\alpha}$.  Note that
$u_1>u_2>u_0$. Whereas the sample size $n_3$ is chosen
to be $M$ in Section~\ref{sec:three}, we now define
\begin{eqnarray*}
\widetilde{n}(\theta)& =& \min\{|\log \alpha | / \inf_{\lambda : u (\lambda) = u_0}
I(\theta, \lambda), |\log \widetilde \alpha | / \inf_{\lambda : u(\lambda) =
  u_2} I(\theta, \lambda)\},\\
n_3&=&n_2\vee\{M'\wedge
\lceil(1+\rho_m)\widetilde{n}(\widehat{\theta}_{n_2})\rceil\},\end{eqnarray*}
where $M<M'\le\widetilde{M}$ and $n_2=m\vee \{M\wedge (1+\rho_m)\widetilde{n}(\widehat{\theta}_m)\}$.  We can regard the test as a group
sequential test with 4 groups and $n_1=m$, $n_4=\widetilde{M}$, but
with adaptively chosen $n_2$ and $n_3$. If the test does not end at the third stage,
continue to the fourth and final stage with sample size
$n_4=\widetilde{M}$. Its rejection and futility
boundaries are similar to those in Section~\ref{sec:three}. Extending
our notation $\Lambda_{i,j}$ in (\ref{eq:1}) to $1\le i\le 4$ and
$0\le j\le 2$, the test stops at stage
$i\le 3$ and rejects $H_0$ if \begin{equation}\label{eq:4} n_i<\widetilde{M},\quad
u(\widehat{\theta}_{n_i})>u_0,\quad\mbox{and}\quad \Lambda_{i,0}\ge b,\end{equation} stops and
accepts $H_0$ if \begin{equation}\label{eq:15}n_i<\widetilde{M},\quad
  u(\widehat{\theta}_{n_i})<u_2,\quad\mbox{and}\quad \Lambda_{i,2}\ge \widetilde{b},\end{equation}
and rejects $H_0$ at stage $i=3$ or $4$ if
\begin{equation}\label{eq:28}n_i=\widetilde{M},\quad u(\widehat{\theta}_{\widetilde{M}})>u_0,\quad\mbox{and}\quad
\Lambda_{i,0}\ge c,\end{equation} accepting $H_0$ otherwise. The thresholds $b, \widetilde{b}$ and $c$ can be defined by
equations similar to (\ref{eq:9})-(\ref{eq:11}) to insure the overall Type I error
probability to be $\alpha$. For example, in place of (\ref{eq:9}),
\begin{equation}\label{eq:16}\sup_{\theta:u(\theta)=u_2} P_\theta \{\mbox{(\ref{eq:15}) occurs for some $i\le
  3$}\}=\widetilde{\eps}\widetilde{\alpha}.\end{equation} The basic
  idea underlying (\ref{eq:16}) is to
  control the Type II error probability at $u_2$ so that the test does not
  lose much power there in comparison with the GLR test
  that has sample size $\widetilde{M}$ (and therefore power
  $1-\widetilde{\alpha}$ at $u_2$).

\subsection{Implementation via Normal Approximation or Monte Carlo}\label{sec:implementation}

\noindent To begin with, suppose the $X_i$ are $N(\theta, 1)$ and
$u(\theta)=\theta$. We write $\theta_j$ instead of $u_j$ and, without 
loss of generality, we shall assume that $\theta_0 = 0$. The
thresholds $b$, $\widetilde b$ and $c$ of the three-stage test in
Section~\ref{sec:three} can be computed by 
solving in succession (\ref{eq:9}), (\ref{eq:10}) and (\ref{eq:11}). Univariate
grid search or Brent's method~\citep{Press92} can be used to
solve each equation.  Since $I(\theta, \lambda) = (\theta - \lambda)^2
/ 2$, we can rewrite (\ref{eq:9}) as $$P_{\theta_1} \{S_m - m
\theta_1 > - (2 \widetilde b m)^{1 \over 2},S_{n_2} - n_2 \theta_1
\leq - (2 \widetilde b n_2)^{1 \over 2}\}
+P_{\theta_1} \{S_m - m \theta_1 \leq - (2 \widetilde b 
m)^{1/2}\}=\widetilde\varepsilon  \widetilde \alpha,$$ and (\ref{eq:10}) and
(\ref{eq:11}) as
\begin{eqnarray*}
P_0 \{S_m / (2m)^{\frac{1}{2}} \geq b^{1 \over 2}\}+ P_0\{\widetilde b^{1 \over 2} < S_m / (2m)^{\frac{1}{2}} < b^{1 \over 2}, 
S_{n_2} / (2 n_2)^{\frac{1}{2}} \geq b^{1 \over 2}, n_2<M\}&=&\varepsilon\alpha,\\
P_0\{\widetilde b^{1 \over 2} < S_m / (2m)^{\frac{1}{2}} < b^{1 \over
  2},n_2<M, \widetilde b^{1 \over 2} < S_{n_2} / (2 n_2)^{\frac{1}{2}} < b^{1 \over
  2}, S_M / (2M)^{\frac{1}{2}} \geq c^{1 \over 2}\}&&\\
+P_0\{\widetilde b^{1 \over 2} < S_m / (2m)^{\frac{1}{2}} < b^{1 \over
  2}, n_2=M,S_M/(2M)^{\frac{1}{2}}\ge c^{1 \over 2}\}&=&(1 - \varepsilon) \alpha.
\end{eqnarray*} The probabilities involving $n_2$ can be computed by
conditioning 
on the value of $S_m/m$, which completely determines the value of
$n_2$, denoted by $k(x)$. For example, the probabilities under
$\theta=0$ can be computed via  
\begin{eqnarray}&&P_0\{S_{n_2} \geq (2 b n_2)^{1 \over 2} | S_m = m x\} = P\{N(0,
1) \geq [2 b k(x) n_2^{1 \over 2} - m x] / [k(x) - m]^{1 \over 2}\},\\
&& P_0 \{S_{n_2} \in dy, S_M \in dz | S_m = mx\} = \varphi_{k(x) -
  m} (y - mx) \varphi_{M - k(x)} (z-y) dydz,\label{eq:27}\end{eqnarray} where $\varphi_v$ is the $N(0, v)$ density function, i.e.,
$\varphi_v(w) = (2 \pi v)^{-{1 \over 2}} \exp(-w^2 / 2v)$. The
probabilities under $\theta_1$ can be computed similarly. Hence
standard recursive numerical integration algorithms can be used to compute the
probabilities in (\ref{eq:9})-(\ref{eq:11}); see
\citet[Section~19.2]{Jennison00}. For the general multiparameter
exponential family, this method can be used to compute the thresholds $b, \widetilde
b$ and $c$ for (\ref{eq:6})-(\ref{eq:8}) since the problem can be
approximated by that of testing a normal mean, as discussed in Section~\ref{sec:three}.

For mid-course modification of the maximum sample
size in Section~\ref{sec:Mtild}, the above recursive numerical algorithm can be modified to handle the randomness of $n_2$ and
$n_3$. The basic idea is that conditional on $S_m/m=x$, the value of
$n_2$ is completely determined as $k(x)$, and conditional on $S_m/m=x$
and $S_{n_2}/n_2=y$, the value of $n_3$ is completely determined as
$h(x,y)$. Therefore, analogous to (\ref{eq:27}), we now have 
\begin{multline}\label{eq:29}
P\{S_{n_3}\in du, S_{\widetilde{M}}\in dw| S_m/m=x, S_{n_2}/n_2=y\}\\
=\varphi_{h(x,y)-k(x)}(u-yk(x))\varphi_{\widetilde{M}-h(x,y)}(w-u)du dw,\end{multline}
and can use bivariate recursive numerical integration. For the general
exponential family, normal approximation to the signed-root likelihood
ratio statistic can again be used.

An alternative to normal approximation is to use Monte Carlo similar to that used in
bootstrap tests. While using Monte Carlo simulations to compute error
probabilities is an obvious idea, it is far from being clear which
distribution from a composite hypothesis should be chosen to simulate
from. Bootstrap theory suggests that we can simulate from the
estimated distribution under the assumed hypothesis as the GLR
statistic is an approximate pivot under that hypothesis. Since the
``estimated distribution'' needs data to arrive at the estimate, we
make use of the first-stage data to determine $b$ and $\widetilde b$,
then we use the second-stage data to determine $c$ for the three-stage
test in Section~\ref{sec:three}.  Specifically, the
Monte Carlo method to determine $b, \widetilde b$ and $c$ proceeds as
follows. At the end of the first stage, compute the maximum likelihood
estimate $\widehat \theta_{m, j}$ under the constraint $u(\theta) =
u_j, j = 0, 1$.  Determine $\widetilde b, b$ and $c$ successively by
solving
\begin{eqnarray}&& P_{\widehat{\theta}_{m,1}}
  \{\mbox{(\ref{eq:7}) occurs for} \ i = 1 \ \mbox{or} \ 2\}=
\widetilde \varepsilon \widetilde \alpha, \label{eq:12}\\
&& P_{\widehat{\theta}_{m,0}} \{\mbox{(\ref{eq:7}) does not occur for}
\ i\le 2, \mbox{and (\ref{eq:6}) occurs for} \ i = 1 \ \mbox{or} \
2\}= \varepsilon \alpha , \label{eq:13}\\
&& P_{\widehat{\theta}_{n_2,0}} \{\mbox{(\ref{eq:6}) and (\ref{eq:7}) do not occur
  for} \ i\le 2, \mbox{and (\ref{eq:8}) occurs}\}=
(1-\varepsilon)\alpha, \label{eq:14}
\end{eqnarray} 
noting that $c$ does not have to be determined until after the second
stage when $n_2$ observations become available for the updated
estimate $\widehat \theta_{n_2, 0}$. The probabilities in
(\ref{eq:12})-(\ref{eq:14}), with $\theta=\widehat{\theta}_{m,1}$ or
$\theta=\widehat{\theta}_{m,0}$ as indicated, can be computed by Monte
Carlo simulations. Similarly, to determine thresholds $b,
\widetilde{b}$ and $c$ of the adaptive test in
Section~\ref{sec:Mtild}, we can use Monte Carlo simulations instead of
normal approximation and numerical integration to compute the
corresponding probabilities.

\section{\large APPLICATIONS AND NUMERICAL EXAMPLES}\label{sec:examples}

\subsection{Application to Single-Arm Phase II Cancer Trials}\label{sec:phaseII}

As pointed out by~\citet[p.~927]{Vickers07}, in a typical phase~II
study of a novel cancer treatment, ``a cohort of patients is treated,
and the outcomes are related to the prespecified target or bar. If the
results meet or exceed the target, the treatment is declared worthy of
further study; otherwise, further development is stopped. This has
been referred to as the `go/ no go' decision. Most often, the outcome
specified is a measure of tumor response, e.g., complete or partial
response using Response Evaluation Criteria in Solid Tumors, expressed
as a proportion of the total number of patients. Response can also be
defined in terms of the proportion who have not progressed or who are
alive at a predetermined time (e.g., one year) after treatment is
started.'' The most widely used designs for these single-arm phase~II
trials are Simon's~(\citeyear{Simon89}) optimal 2-stage designs, which
allow early stopping of the trial if the treatment has not shown
beneficial effect that is measured by a Bernoulli proportion. These
designs are optimal in the sense of minimizing the expected sample
size under the null hypothesis of no viable treatment effect, subject
to Type~I and II error probability bounds. Given a maximum sample
sample size $M$, Simon considered the design that stops for futility
after $m<M$ patients if the number of patients exhibiting positive
treatment effect is $r_1 (\le m)$ or fewer, and otherwise treats an
additional $M-m$ patients and rejects the treatment if and only if the
number of patients exhibiting positive treatment effect is $r_2(\le
M)$ or fewer. Simon's designs require that a null proportion $p_0$,
representing some ``uninteresting'' level of positive treatment
effect, and an alternative $p_1>p_0$ be specified.  The null
hypothesis is $H_0:p\le p_0$, where $p$ denotes the probability of
positive treatment effect. The Type~I and II error probabilities
$P_{p_0}\{\mbox{Reject $H_0$}\}$, $P_{p_1}\{\mbox{Accept
  $H_0$}\}$ and the expected sample size $E_{p_0}N$ can be
computed for any design of this form, which can be
represented by the parameter vector $(m,M,r_1,r_2)$. Using computer
search over these integer-valued parameters, \citet{Simon89} tabulated the
optimal designs in his Tables~1 and 2 for different values of
$(p_0,p_1)$.

Whether the new treatment is declared promising in a phase~II trial
depends strongly on the prescribed $p_0$ and $p_1$. In their
systematic review of $134$ papers reporting phase~II trials in
\textit{J.\ Clin.\ Oncology}, \citet{Vickers07} found 70 papers
referring to historical data for their choice of the null or
alternative response rate, and that nearly half (i.e., 32) of these
papers did not cite the source of the historical data used, while only
9 gave clearly a single historical estimate of their choice of $p_0$.
Moreover, no study ``incorporated any statistical method to account
for the possibility of sampling error or for differences in case mix
between the phase~II sample and the historical cohort.'' The adaptive
designs in Section~\ref{sec:three} applied to this setting do not
require the specification of the alternative $p_1$, a desirable
property to prevent well-intentioned but misguided practitioners from
choosing $p_1$ artificially small to inflate the appearance of a
positive treatment effect, if one exists; uncertainty in the choice of
$p_0$ is also an important issue and is addressed in
Section~\ref{sec:phaseII-III}. For now, assume $p_0$ to be given along
with initial and maximum sample sizes $m$ and $M$.  The adaptive test
takes $p_1$ to be the alternative where the FSS test, with Type~I
error probability $\alpha$ at $p_0$, has power $1-\beta$, i.e., the
solution of $F_{M,p_1}(F^{-1}_{M,p_0}(1-\alpha))=\beta$, where
$F_{M,p}$ is the distribution function of the $\mbox{Bin}(M,p)$
distribution. The GLR statistic~(\ref{eq:1}) at the $i$th stage is
$$n_i
\left[\widehat{p}_{n_i}\log\left(\frac{\widehat{p}_{n_i}}{p_j}\right)+(1-\widehat{p}_{n_i})\log\left(\frac{1-\widehat{p}_{n_i}}{1-p_j}\right)\right]$$
for $j=0$ or $1$, and the critical values $b, \widetilde{b}$ and $c$
are chosen to satisfy (\ref{eq:9})-(\ref{eq:11}). Because of
discreteness of the binomial distribution it may be impossible to
satisfy (\ref{eq:9})-(\ref{eq:11}) exactly, in which case (\ref{eq:9})
is satisfied approximately and (\ref{eq:10})-(\ref{eq:11}) are
satisfied conservatively. The stopping rule defined by
(\ref{eq:6})-(\ref{eq:8}) may alternatively be stated in terms of the
number of cumulative successes $S_{n_i}$ at the $i$th stage. Table~1
describes the adaptive design (denoted by ADAPT) and
Simon's~(\citeyear{Simon89}) optimal 2-stage design (denoted by Sim2)
for two choices of $m,M,\alpha,\beta,p_0$ and Table~2 contains their
operating characteristics, computed exactly using the
$\mbox{Bin}(n,p)$ distribution. ADAPT has expected sample size close
to Sim2 for $p$ near $p_0$, and smaller sample size when $p$ is
roughly midway between $p_0$ and $p_1$ or is larger; $p_1=.3$ in the
top panel of Table~2 and $p_1=.44$ in the bottom panel. It is not
surprising that the expected sample size of Sim2 increases with $p$
since Sim2 only stops early for futility. The expected number of
stages shows a similar pattern, while their power functions are nearly
identical. Note that even though ADAPT has a maximum of three stages,
its expected number of stages is less than 2 for all $p$ and usually
close to 1.

\begin{table}[t!]
\noindent\textbf{Table 1.} Description of ADAPT and Sim2 for two cases
\begin{center}
\begin{tabular}{r@{\hspace{1cm}}ll}
\hline
$S_m$&\hspace{2.5cm}ADAPT&Sim2 ($p_1=.3$ or $.44$)\\\hline
\multicolumn{3}{c}{(a) $m=10$, $M=29$, $p_0=.1$, $\alpha=.05$,
  $\beta=.2$}\\
$\le 1$&Accept $H_0.$&Accept $H_0$.\\
2&$n_2=M$; reject $H_0$ if $S_{n_2}\ge 6$.&$n_2=M$ and\\
3&$n_2=20$&\hspace{.5cm}reject $H_0$ if $S_M\ge
6$.\\
&(i) If $S_{n_2}\le 3$, accept $H_0$.&\\
&(ii) If $S_{n_2}\ge 6$, reject $H_0$.&\\
&(iii) If $4\le S_{n_2}\le 5$ and $S_M\ge 6$, reject $H_0$.&\\
$\ge 4$&Reject $H_0$.&$n_2=M$; rej.\ $H_0$ if $S_M\ge 6$.\\
\\
\multicolumn{3}{c}{(b) $m=30$, $M=82$, $p_0=.1$, $\alpha=\beta=.1$}\\
$\le 8$&Accept $H_0$.&Accept $H_0$.\\
9&$n_2=57$&Accept $H_0$.\\
&(i) If $S_{n_2}\le 19$, accept $H_0$.&\\
&(ii) If $S_{n_2}\ge 24$, reject $H_0$.&\\
&(iii) If $20\le S_{n_2}\le 23$ and $S_M\ge 32$, reject $H_0$.&\\
$10-13$&$n_2=M$; reject $H_0$ if $S_{n_2}\ge 31$.&$n_2=M$ and\\
$\ge 14$&Reject $H_0$.&\hspace{.5cm}reject $H_0$ if $S_M\ge
30$.\\\hline
\end{tabular}
\end{center}
\end{table}

\begin{table}
\noindent\textbf{Table 2.} Expected sample size, power (in parentheses) and expected number of
stages (in brackets) of phase~II designs
\begin{center}
 \begin{tabular}{ccccccccc}
\hline
$p$&&\multicolumn{3}{c}{ADAPT}&&\multicolumn{3}{c}{Sim2}\\\hline
&&\multicolumn{7}{c}{(a) $m=10$, $M=29$, $p_0=.1$, $\alpha=.05$, $\beta=.2$}\\
.05&&11.6&(.3\%)&[1.1]&&11.6&(.2\%)&[1.1]\\
$p_0=.1$&&14.5&(5\%)&[1.3]&&15.0&(4.7\%)&[1.3]\\
.2&&18.8&(43.3\%)&[1.6]&&21.9&(43.1\%)&[1.6]\\
$p_1=.3$&&18.1&(79.4\%)&[1.6]&&26.1&(79.6\%)&[1.8]\\
.4&&14.8&(94.9\%)&[1.4]&&28.1&(95.0\%)&[2.0]\\
.5&&12.1&(98.9\%)&[1.2]&&28.8&(98.9\%)&[2.0]\\
.6&&10.1&(99.9\%)&[1.0]&&29.0&(99.9\%)&[2.0]\\
&&\multicolumn{7}{c}{(b) $m=30$, $M=82$, $p_0=.3$, $\alpha=\beta=.1$}\\
.2&&34.9&(.3\%)&[1.1]&&33.2&(.03\%)&[1.1]\\
$p_0=.3$&&51.8&(10.0\%)&[1.5]&&51.4&(10.0\%)&[1.4]\\
.35&&60.4&(35.0\%)&[1.7]&&63.4&(36.2\%)&[1.6]\\
$p_1=.44$&&52.9&(88.7\%)&[1.5]&&77.7&(87.8\%)&[1.9]\\
.5&&42.4&(98.4\%)&[1.3]&&80.9&(97.5\%)&[2.0]\\
.6&&31.9&(99.9\%)&[1.0]&&82.0&(99.9\%)&[2.0]\\\hline
\end{tabular}
\end{center}
\end{table}

\subsection{Extension to Randomized Phase II Cancer Trials}\label{sec:phaseII-III}

As noted by~\citet{Vickers07}, uncertainty in the choice of $p_0$ and
$p_1$ can increase the likelihood that (a) a treatment with no viable
positive treatment effect proceeds to phase~III, or (b) a treatment
with positive treatment effect is abandoned at phase~II. To circumvent
the problem of choosing an artificially small $p_0$, either
intentionally by a practitioner wanting to give the treatment the
``best chance'' of showing a positive effect if one exists, or
unintentionally because of inaccurate information about the control,
\citet{Ellenberg85} proposed to perform a controlled 2-arm phase~II
trial in which patients are randomized into both treatment and control
groups. After randomizing $2n_1$ patients into treatment and control
arms, the trial is stopped for futility when the number in the
treatment group showing positive effect is not greater than in the
control group.  Otherwise, the trial continues until a total of $2n_2$
patients have been randomized into the study and then a standard fixed
sample binomial test is performed. Letting $p$ and $q$ denote the
probability of positive effect in the treatment and control groups,
respectively, \citet{Thall88} subsequently chose $n_1, n_2$ and two
other design parameters $y_1$ and $y_2$, described below, to minimize the ``average'' expected sample size
\begin{equation}\label{eq:25}
  \mbox{AvSS}=\frac{1}{2}[E(N|p=q)+E(N|p=q+\delta)]\end{equation} subject to Type
  I and II error probability constraints $\alpha$ and $\beta$. The
  two-stage test of $H_0: p\le q$ stops for futility after the first stage if an approximately
  normally distributed test statistic $Z_1$, based on the
  2-population binomial data, is no greater than $y_1$, otherwise
  continuing with a second stage and rejecting $H_0$ if
  $Z_2>y_2$. Because the
  expectations in (\ref{eq:25}) depend on $p$ and $q$, both must be
  specified to calculate the
  trial design.

The adaptive designs of Section~\ref{sec:design} apply naturally to
this 2-arm setting. The GLR statistic~(\ref{eq:1}) at the
$i$th stage is
$$n_i\left[
  \widehat{p}_i\log\left(\frac{\widehat{p}_i}{\widehat{p}_{\delta_j}}\right)
  -(1-\widehat{p}_i)\log\left(\frac{1-\widehat{p}_i}{1-\widehat{p}_{\delta_j}}\right)
  +\widehat{q}_i\log\left(\frac{\widehat{q}_i}{\widehat{p}_{\delta_j}-\delta_j}\right)
  -(1-\widehat{q}_i)
  \log\left(\frac{1-\widehat{q}_i}{1-\widehat{p}_{\delta_j}+\delta_j}\right)\right],$$
where $\widehat{p}_i, \widehat{q}_i$ are the fraction of successes in
the treatment and control arms and
$\widehat{p}_{\delta_j}$ is the MLE of $p$ under $p-q=\delta_j$, $j=0,
  1$, where
$\delta_0=0$ and $\delta_1$ is the implied alternative (see the example below). The boundaries $b,\widetilde{b}$ and $c$
  can be computed using a normal approximation or Monte Carlo
  simulations as
  described in Section~\ref{sec:implementation}; the latter, with
  1 million simulations for each probability calculation, is used in
  the following comparative study. 

  \citet{Thall94} describe a phase II trial of fludarabine + ara-C +
  granulocyte colony stimulating factor~(G-CSF) for treatment of acute
  myelogenous leukemia~(AML). Although this trial was designed using
  other methods described in their paper, we use it here as a real
  setting to compare the adaptive design with Thall et
  al.'s~(\citeyear{Thall88}) design. The standard therapy for AML at
  the time was fludarabine + ara-C alone, and in a preceding study, 22
  out of 45 patients achieved complete remission of the leukemia, the
  clinical endpoint of interest, suggesting an initial estimate of $q$
  of .5. The study was conducted to detect an increase in remission
  rate of $\delta=.2$. For $\alpha=.05$ and $\beta=.2$, Thall et
  al.'s~(\citeyear{Thall88}) optimal 2-stage design (denoted by Opt2)
  for detecting a 20\% improvement when the control remission rate $q$
  is .5 has first stage of 33 per arm, followed by a second stage of
  45 per arm if $Z_1>y_1=.356$, and rejecting the null hypothesis
  after the second stage if $Z_2>y_2=1.584$. The adaptive design
  (denoted by ADAPT) with first stage $m=25$ and maximum sample 78 per
  arm, the same as Opt2, uses boundaries $b=2.12, \widetilde{b}=1.03$
  and $c=1.56$ for $\alpha=.05, \widetilde{\alpha}=.2$ and
  $\eps=\widetilde{\eps}=1/2$. Table~3 contains the operating
  characteristics of ADAPT and Opt2 for a variety of treatment and
  control response rates $(p,q)$ around $(.7,.5)$. Each result is
  based on 100,000 simulations. ADAPT has substantially smaller
  expected sample size than Opt2. This is in part because Opt2 only
  stops early for futility, although the parameters of Opt2 in this
  case are chosen to minimize (\ref{eq:25}), yet there is substantial
  savings both when $p=q$ and $p=q+\delta$. ADAPT and Opt2 have
  similar expected number of stages near the null hypothesis, with
  ADAPT decreasing as $p-q$ increases while Opt2 steadily increases to
  2, again due to its early stopping only for futility. The power
  functions of the tests are similar, with Opt2 having slightly higher
  power.  Note that the Type I error probability of Opt2 is inflated
  above $\alpha=.05$ at $p=q=.5$ due to the normal approximations to
  $Z_i$ used to compute the design parameters.

\begin{table}[t!]
\noindent\textbf{Table 3.}  Expected sample size, power (in parentheses), expected number of
  stages (in brackets) and average expected sample size (\ref{eq:25})
  (denoted by AvSS) of 2-arm phase~II designs
\begin{center}
 \begin{tabular}{cccccccccc}
\hline
$q$&$p$&&\multicolumn{3}{c}{ADAPT}&&\multicolumn{3}{c}{Opt2}\\\hline
.4&.3&&33.3&(0.4\%)&[1.1]&&37.8&(0.2\%)&[1.1]\\
&.4&&46.1&(5.3\%)&[1.5]&&48.9&(5.3\%)&[1.4]\\
&.5&&57.5&(32.3\%)&[1.8]&&63.3&(35.6\%)&[1.7]\\
&.6&&56.4&(76.0\%)&[1.8]&&73.5&(78.9\%)&[1.9]\\
&.7&&43.8&(97.0\%)&[1.5]&&77.3&(97.7\%)&[2.0]\\
\multicolumn{2}{c}{AvSS}&&51.3&&&&61.2&&\\
&&&&&&&&&\\

.5&.4&&34.7&(0.4\%)&[1.2]&&38.2&(0.2\%)&[1.1]\\
&.5&&47.3&(5.0\%)&[1.5]&&49.0&(5.6\%)&[1.4]\\
&.6&&57.5&(32.2\%)&[1.8]&&63.3&(35.5\%)&[1.7]\\
&.7&&55.1&(77.8\%)&[1.8]&&73.7&(80.4\%)&[1.9]\\
&.8&&41.0&(97.6\%)&[1.4]&&77.5&(98.2\%)&[2.0]\\
\multicolumn{2}{c}{AvSS}&&51.2&&&&61.4&&\\
&&&&&&&&&\\

.6&.5&&34.7&(0.4\%)&[1.2]&&38.2&(0.2\%)&[1.1]\\
&.6&&46.0&(5.2\%)&[1.5]&&48.9&(5.3\%)&[1.4]\\
&.7&&55.8&(33.2\%)&[1.7]&&63.3&(35.6\%)&[1.7]\\
&.8&&52.3&(81.1\%)&[1.7]&&74.4&(84.2\%)&[1.9]\\
&.9&&35.9&(98.5\%)&[1.3]&&77.8&(99.4\%)&[2.0]\\
\multicolumn{2}{c}{AvSS}&&49.2&&&&61.7&&\\\hline
\end{tabular}
\end{center}
\end{table}

\subsection{Comparison with Adaptive Tests for Difference of Means with Unknown Variances}\label{sec:var_unknown}

Let $X_1,X_2,\ldots$ and $Y_1,Y_2,\ldots$ be independent normal
observations with unknown means $\mu_X, \mu_Y$ and variances
$\sigma_X^2, \sigma_Y^2$, respectively. Even if the variances
$\sigma_X^2, \sigma_Y^2$ are assumed equal, no fixed sample size test
of the hypothesis $\mu_X\le\mu_Y$ can achieve specified error
probabilities $\alpha$, $\tilde{\alpha}$ at $\mu_X-\mu_Y$ equal to 0
and some specified $\delta>0$ without knowing the true value of the
variance, as demonstrated by~\citet{Dantzig40}. To overcome this
difficulty in the case of a single normal mean, \citet{Stein45}
proposed a two-stage procedure that uses the first-stage sample to
estimate the variance.  The total sample size is then determined as a
function of this estimate so that the test statistic, which uses the
data from both stages to estimate the mean but only the first stage to
estimate the variance, is exactly $t$-distributed under
$\mu_X-\mu_Y=0$ at which the test has Type I error $\alpha$; the power
of the test at $\mu_X-\mu_Y=\delta$ is strictly larger than
$1-\widetilde{\alpha}$. Stein's procedure can be easily extended to
the mean difference problem if $\sigma_X^2=\sigma_Y^2$, and has been
modified for clinical trials by~\citet{Wittes90}, \citet{Birkett94},
and~\citet{Denne99}. These modifications of Stein's procedure
use the overall sample to estimate both the mean difference and the
common variance in modifying Stein's test statistic, and therefore may
not maintain the Type I error probability at the prescribed level.
\citet{Wittes90} also assume a prior estimate $\sigma_0^2$ (e.g., from
previous studies in the literature) of the common variance to modify
Stein's formula for the total sample size.

In Table~4 we compare the performance of the adaptive test in
Section~\ref{sec:three} (denoted by ADAPT) with Stein's test, denoted
by S, and
the modified versions of~\citet[denoted by WB]{Wittes90},
\citet[denoted by BD]{Birkett94} and~\citet[denoted by DJ]{Denne99} in
the context of a phase II hypercholesterolaemia treatment efficacy
trial described by~\citet{Facey92}.  In this trial, patients were
randomized into treatment and placebo groups and serum cholesterol
level reductions, $X_i$ and $Y_i$, assumed to be normally distributed,
were measured after four weeks of treatment.  A difference in
reductions of serum cholesterol levels, in mmol/liter, between the
treatment and placebo groups of $1.2$ was of clinical interest. Based
on previous studies, it was anticipated that the standard deviation of
the reductions would be about 0.7 for both groups.  If the standard
deviation were known to be $\sigma_0=0.7$, the size of the fixed
sample $t$-test with error probabilities $\alpha=\tilde{\alpha}= .05$
at mean difference 0 and $\delta=1.2$ is 9 per group. Following
\citet{Denne99}, we assume a first-stage per-group sample size of
$m=5$, being approximately half of 9.  If the standard deviation were
in fact $2\sigma_0 =1.4$, the per-group sample size of the same
$t$-test is 31, which we take as a reasonable maximum sample size $M$
for our three-stage test with $\rho_m=.1$ and
$\varepsilon=\widetilde{\varepsilon}=1/3$. Table~4 contains the power
and per-group expected sample size of ADAPT and the aforementioned
procedures in the literature, evaluated by 100,000 simulations at
various values of $\mu_X-\mu_Y\in[0,\delta]$ and
$\sigma=\sigma_X=\sigma_Y$. Whereas the Stein-type tests require this
assumption of equal variances, the three-stage tests defined in
Section~\ref{sec:three} do not, so for comparison we also include in
Table~4 the three-stage test that does not assume $\sigma_X=\sigma_Y$,
which we denote by ADAPT$_{\ne}$. The Kullback-Leibler information
number $I=\min\{ I((\mu_X,\mu_Y,\sigma^2), (\widetilde{\mu}_X,
\widetilde{\mu}_Y,\widetilde{\sigma}^2)):
\widetilde{\mu}_X-\widetilde{\mu}_Y=0\}$, where
$I(\theta,\widetilde{\theta})$ is defined in (\ref{eq:2}), is also
reported in the first column. When the true standard deviations
$\sigma_X$ and $\sigma_Y$ are equal to the specified value $\sigma_0$,
ADAPT and ADAPT$_{\ne}$ have similar power but smaller expected sample
size than the other tests for values of $\mu_X-\mu_Y$ near 0 and
$\delta$.  When the standard deviations $\sigma_X$ and $\sigma_Y$ are
larger than the specified value $\sigma_0$, the adaptive tests have
much smaller expected sample sizes than the Stein-type tests, whose
second-stage sample size increases without bound as a function of the
first-stage sample variance; in particular, see the last three rows of
Table~4.

\begin{table}[th!]
\noindent\textbf{Table 4.}  Power and per-group expected sample size of tests of $H_0:\mu_X\le\mu_Y$
\begin{center}
  \begin{tabular}{ccccccc}
\hline
$(\mu_X-\mu_Y,\sigma)$&&&&&\\
$I$&\raisebox{1.5ex}[0pt]{S}&\raisebox{1.5ex}[0pt]{WB}&\raisebox{1.5ex}[0pt]{BD}&\raisebox{1.5ex}[0pt]{DJ}&\raisebox{1.5ex}[0pt]{ADAPT}&\raisebox{1.5ex}[0pt]{ADAPT$_{\ne}$}\\\hline
$(0,\sigma_0/2)$&5.0\%&4.9\%&5.0\%&5.0\%&1.6\%&1.8\%\\
$I=0$&5.0&9.0&5.0& 5.0&5.0&5.0\\
\\
$(0,\sigma_0)$&5.0\%&5.4\%&6.0\%&5.6\%&4.0\%&4.1\%\\
$I=0$&10.1&10.1&10.1&10.3&9.4&10.2\\
\\
$(\delta/2,\sigma_0)$&53.0\%&59.5\%&55.4\%&57.3\%&65.0\%&68.1\%\\
$I=.169$&10.1&10.1&10.1&10.3&15.5&13.7\\
\\
$(\delta,\sigma_0)$&96.6\%&98.0\%&95.8\%&96.4\%&97.8\%&98.5\%\\
$I=.551$&10.1&10.1&10.1&10.3&9.4&8.0\\
\\
$(0,2\sigma_0)$&5.0\%&5.5\%&5.5\%&4.6\%&5.0\%&5.3\%\\
$I=0$&38.2&38.2&38.2&30.7&22.1&22.7\\
\\
$(\delta/2,2\sigma_0)$&50.3\%&50.0\%&49.7\%&53.8\%&44.0\%&44.5\%\\
$I=.045$&38.1&38.2&38.2&30.7&25.5&26.0\\
\\
$(\delta,2\sigma_0)$&95.2\%&89.1\%&91.3\%&92.9\%&91.9\%&91.4\%\\
$I=.169$&38.1&38.2&38.2&30.7&22.1&20.2\\
\\
$(0,3\sigma_0)$&5.0\%&5.3\%&5.3\%&4.6\%&5.1\%&5.2\%\\
$I=0$&85.2&85.2&85.3&67.7&26.3&26.7\\
\\
$(0,5\sigma_0)$&5.0\%&5.4\%&5.4\%&4.7\%&5.2\%&5.1\%\\
$I=0$&236&235&236&186&27.6&28.6\\
\\
$(0,10\sigma_0)$&5.0\%&5.4\%&3.8\%&4.9\%&5.1\%&5.1\%\\
$I=0$&942&940&943&754&28.7&29.3\\\hline
\end{tabular}\end{center}\end{table}

An alternative approach to Stein-type designs has been used by
\citet{Proschan95} and~\citet{Li02}, who simply replace the $\sigma^2$
in their two-stage tests that assume known variance with its current
estimate at each stage.  To compare ADAPT with these tests, which rely
on stable variance estimates, we allow a larger first-stage sample
size of $m=20$.  Table~5 contains the power and per-group expected
sample size of Proschan \& Hunsberger's test (denoted by PH), two
choices of the early stopping boundaries $(h,k)$ in Table~1 of
\citet{Li02} for their test, which we denote by L1 and L2, and our
three-stage test ADAPT, for various values of $\mu_X-\mu_Y$ and
$\sigma$, each entry being the result of 100,000 replications.  To
compare these tests on equal footing we have chosen the maximum sample
size $M=121$ for ADAPT because this is the maximum sample size of L1
and is quite close to the maximum sample sizes of PH and L2, which are
122 and 104, respectively. The PH, L1 and L2 tests are designed to
achieve Type~I error probability $.05$ and they choose the sample size of
their second stage based on a conditional power level of 80\%.  The
threshold values $b=2.68, \widetilde{b}=1.75, c=1.75$ used by ADAPT
are thus computed using $\alpha=.05,\tilde{\alpha}=.20$. The results
in Table~5 show that the true power of L1, L2 and PH falls well below
their nominal conditional power level of 80\%. When $\sigma=2$, the
L1, L2 and PH tests have power less than 50\% for all values of
$\mu_X-\mu_Y$ considered, which is caused by stopping prematurely for
futility at the end of the first stage; see in particular the rows in
Table~5 that correspond to $\mu_X-\mu_Y=0$.  Since the conditional
power criterion is not valid when the estimated difference of means is
near zero, the L and PH tests must stop for futility when this occurs
even though the true difference of means may be substantially greater
than zero.

\begin{table}[t!]
\noindent\textbf{Table 5.} Maximum sample size $M$, power and
per-group expected sample size for the tests L1 and L2 of Li et al.,
Proschan \& Hunsberger (PH), and ADAPT
\begin{center}
  \begin{tabular}{ccccc}
\hline
$(\mu_X-\mu_Y,\sigma)$&&&&\\
$I$&\raisebox{1.5ex}[0pt]{L1}&\raisebox{1.5ex}[0pt]{L2}&\raisebox{1.5ex}[0pt]{PH}&\raisebox{1.5ex}[0pt]{ADAPT}\\\hline
$(0,1)$&5.5\%&5.3\%&5.5\%&4.8\%\\
$I=0$&26.3&25.5&25.9&56.5\\
\\
$(0,2)$&5.3\%&5.3\%&5.4\%&5.4\%\\
$I=0$&26.2&26.3&25.9&93.5\\
\\
$(1/4,1)$&29.9\%&29.3\%&29.0\%&48.3\%\\
$I=.016$&32.8&31.0&31.6&76.1\\
\\
$(3/8,1)$&49.5\%&48.7\%&48.3\%&77.5\%\\
$I=.035$&34.5&32.7&33.1&73.5\\
\\
$(1/2,1)$&67.8\%&66.4\%&66.4\%&92.8\%\\
$I=.061$&34.3&32.8&32.7&63.3\\
\\
$(1/2,2)$&12.0\%&29.9\%&28.9\%&56.1\%\\
$I=.016$&29.1&32.9&31.7&98.7\\
\\
$(3/4,2)$&49.8\%&48.5\%&48.1\%&85.6\%\\
$I=.035$&34.5&32.7&33.2&87.0\\\hline
\end{tabular}\end{center}
\end{table}

\subsection{Comparison of Tests Allowing Mid-course Modification of
  Maximum Sample Size}\label{sec:Mtild_eg}

As pointed out in Section~\ref{sec:intro}, \citet{Cui99} have proposed a
method to modify the group size of a given group sequential test of
$H_0:\theta\le0$ in response to protocol amendments during interim analyses.  In the example
considered by~\citet[p.~854]{Cui99}, the maximum sample
size is initially $M=125$ for detecting $\theta_1=.29$ with power .9
and $\alpha=.025$, but can be subsequently increased up to
$\widetilde{M}=500$; their sample sizes are twice as large because
they consider variance 2. They consider modifying the group size at the end of a
given stage $L$ if the ratio of conditional power at the observed
alternative $\widehat{\theta}_{n_L}$ to the conditional power at
$\theta_1$ is greater than 1 or less than .8, in which case the group
size is then modified so that the new maximum sample size is
\begin{equation}
\widetilde{M}\wedge M(\theta_1/\widehat{\theta}_{n_L})^2.  \label{eq:26}
\end{equation} 

If (\ref{eq:26}) is
less than the already sampled $n_L$, error spending can be used to end
the trial.  The crux of this method is that the original critical
values can be used for the weighted test statistic without changing
the Type I error probability regardless of how the sample size is
changed. Table~6 compares their proposed adaptive group sequential
tests with FSS tests, standard (non-adaptive) group sequential tests,
and the adaptive test described in Section~\ref{sec:Mtild}. Each
result is based on 100,000 simulations. All
adaptive tests in Table~6 use the first-stage sample size $m=25$,
maximum sample size initially $M=125$ with the possibility of
extension up to $\widetilde{M}=500$, and Type I error probability not
exceeding $\alpha=.025$, matching the setting considered in Section~2
of~\citet{Cui99}. Since the maximum sample size can vary between
$M=125$ and $\widetilde{M}=500$, the two relevant implied alternatives
are $\theta_1=.29$, where FSS$_{125}$ has power
$1-\widetilde{\alpha}=.9$, and $\theta_2=.15$, where FSS$_{500}$ has
power .9. The values of the user-specified parameters of the tests in
Table~6 are summarized below.

\newpage
\begin{center}
INSERT TABLE 6 ON THIS PAGE 
\end{center}
\newpage

\begin{itemize}
\item ADAPT: The adaptive test, described
in Section~\ref{sec:Mtild}, that uses $b=3.48$,
$\widetilde{b}=2.1$ and $c=2.31$ corresponding to
$\varepsilon=\widetilde{\varepsilon}=1/2$, $\rho_m=.1$ and $M'=250$.
\item FSS$_{125}$, FSS$_{500}$: The FSS tests having sample sizes
  $M=125$ and $\widetilde{M}=500$, respectively.
\item OBF$_{SC}^5$: A one-sided O'Brien-Fleming
group sequential test having five groups of size 100 and that uses stochastic
curtailment futility stopping ($\gamma=.9$ in Section~10.2 of~\citet{Jennison00}) with reference alternative
$\theta_2=.15$; see below.
\item C$^4$, C$^5$: Two versions the  adaptive group
sequential test of~\citet{Cui99} that adjusts the group size at the end of the first
stage; C$^4$ uses four stages and C$^5$ uses five stages.
\item C$_{SC}^5$, C$_{PF}^5$: Two modifications of C$^5$ to allow for
  futility stopping; C$_{SC}^5$ uses stochastic curtailment futility
stopping ($\gamma=.9$ in Section~10.2 of~\citet{Jennison00}) and C$_{PF}^5$ uses power family
  futility stopping ($\Delta=1$ in Section~4.2 of~\citet{Jennison00}). Both C$_{SC}^5$ and C$_{PF}^5$ use
  reference alternative 
  $\theta_2=.15$.
\end{itemize}
Since OBF$_{SC}^5$,
C$_{SC}^5$ and C$_{PF}^5$ have maximum sample size
$\widetilde{M}=500$, the futility stopping
boundaries of these tests are designed to
have power .9 at $\theta_2$. We have also included
C$^4$ because our adaptive test uses no
more than four stages.  The tests are evaluated at the $\theta$ values where
FSS$_{125}$ has power .01, .025, .7, .8 and .9, and where FSS$_{500}$
has power .7, .8 and .9.

Even though the C tests have maximum sample size $\widetilde{M}=500$,
they are underpowered at $0<\theta\le\theta_2$, the alternative
implied by $\widetilde{M}$, when compared with ADAPT, FSS$_{500}$ and
OBF$_{SC}^5$. In particular, the C tests have power less than .65 at
$\theta_2$.  Since C$^4$ and C$^5$ use no futility stopping, this
suggests that their updated maximum sample size (\ref{eq:26}) (with
$L=1$) has contributed to the power loss. The large expected sample
sizes of C$^4$ and C$^5$ at $\theta\le0$ reveal another problem with
this sample size updating rule: It does not consider the sign of
$\widehat{\theta}_{m}$; a negative value of $\widehat{\theta}_{m}$
could result in the same sample size modification as a positive one,
causing a large increase in the group size when it should be decreased
toward futility stopping. ADAPT has only a slight loss of power in
comparison with FSS$_{500}$ and the five-stage OBF$_{SC}^5$ at
$\theta>0$, with substantially smaller expected sample size.  The mean
number of stages of ADAPT at $\theta_1=.29$ shows that it behaves like
a two- or three-stage test there. OBF$_{SC}^5$, on the other hand, has
the largest expected sample size at $\theta\ge0$ of the tests in
Table~6 other than FSS$_{125}$.

\section{\large ASYMPTOTIC THEORY AND A MODIFIED CONDITIONAL POWER TEST}\label{sec:asymptotics}

In this section we prove the asymptotic optimality of the adaptive
tests in Sections~\ref{sec:three} and \ref{sec:Mtild}. The proof also
sheds light on how the two-stage conditional power tests, which are
shown to be severely under-powered in Section~\ref{sec:examples}, can
substantially increase their power by adding a third stage. These
modified conditional power tests are still not asymptotically
efficient because they try to mimic the optimal FSS test when the
alternative is given whereas the adaptive test of
Section~\ref{sec:three} tries to mimic the SPRT instead, assuming
$H_0$ to be simple.  

\bigskip

\noindent\textbf{Theorem~4.1.}  \emph{Let $N$ denote the sample size of
the three-stage GLR test in Section~\ref{sec:three}, with $m$, $M$ and $m\vee\{M\wedge\lceil (1 +
\rho_m) n(\widehat \theta_m)\rceil\}$ being the possible values of
$N$. Let $T$ be the sample size of any test of $H_0: u(\theta)\le u_0$
versus $H_1: u(\theta)\ge u_1$,
sequential or otherwise, which takes at least $m$ and at most $M$
observations and whose Type I and Type II error probabilities do not
exceed $\alpha$ and $\widetilde{\alpha}$, respectively. Assume that
$\log \alpha \sim \log \widetilde \alpha$,
\begin{equation}\label{eq:17}m/|\log \alpha |
\rightarrow a, \ M / |\log \alpha | \rightarrow A, \rho_m \rightarrow
0\;\; \mbox{but}\;\; m^{\frac{1}{2}}\rho_m/(\log m)^{1/2}\rightarrow
\infty \end{equation} as $\alpha + \widetilde \alpha \rightarrow 0$,
with $0 < a < A$. Then for every fixed $\theta$, as $\alpha +
\widetilde \alpha \rightarrow 0$, \begin{eqnarray}E_\theta (N) &\sim&
m \vee (M \wedge |\log \alpha | / \{ \inf_{\lambda : u(\lambda) = u_0}
I(\theta, \lambda) \vee \inf_{\lambda : u(\lambda) = u_1} I(\theta,
\lambda)\}),\label{eq:18}\\ E_{\theta}(T)&\ge&
(1+o(1))E_{\theta}(N).\end{eqnarray}}

\noindent\textit{Proof.} Let $\Theta_0=\{\theta:u(\theta)\le u_0\}$,
$\Theta_1=\{\theta:u(\theta)\ge u_1\}$. By (\ref{eq:30}), for $i=0,1$,
\begin{equation}\label{eq:31}\inf_{\lambda\in\Theta_i}
  I(\theta,\lambda)=I_i(\theta),\quad\mbox{where}\quad
  I_i(\theta)=\inf_{\lambda:u(\lambda)=u_i}
  I(\theta,\lambda).\end{equation} Take any $\lambda\in\Theta_0$ and
  $\widetilde{\lambda}\in\Theta_1$. In view of (\ref{eq:31}) and
  Hoeffding's~(\citeyear{Hoeffding60}) lower bound for the expected
  sample size $E_\theta (T)$ of a test that has error probabilities
  $\alpha$ and $\widetilde{\alpha}$ at $\lambda$ and
  $\widetilde{\lambda}$ and takes at least $m$ and at most $M$
  observations, 
\begin{equation}\label{eq:32}E_\theta (T)\ge m\vee\left\{M\wedge
  \frac{(1+o(1)) |\log\alpha|}{I_0(\theta)\vee
  I_1(\theta)}\right\}\end{equation} as $\alpha +\widetilde{\alpha}\To
  0$ such that $\log \alpha\sim \log \widetilde{\alpha}$. We next show
  that the asymptotic lower bound in (\ref{eq:32}) is attained by
  $N$. Since $m\sim a|\log \alpha|$ and $M\sim A|\log \alpha|$ and
  since the thresholds $b, \widetilde{b}$ and $c$  are defined by
  (\ref{eq:9})-(\ref{eq:11}), we can use an argument similar to the
  proof of Theorem~2(ii) of~\citet[p.~525]{Lai04} to show that
  (\ref{eq:18}) holds. In fact, the second-stage sample size of the
  adaptive test is a slight inflation of the Hoeffding-type lower
  bound (\ref{eq:32}) with $\theta$ replaced by the maximum likelihood
  estimate $\widehat{\theta}_m$ at the end of the first stage. The
  assumption $\rho_m\To 0$ but $\rho_m \succ m^{-1/2}(\log m)^{1/2}$
  is used to accommodate the difference between $\theta$ and its
  substitute $\widehat{\theta}_m$, noting that as $m\To\infty$,
  $$P_\theta \{\sqrt{m}|\widehat{\theta}_m-\theta|\ge
  r(\log m)^{1/2}\}=o(m^{-1})$$ is $r$ if sufficiently large, by
  standard exponential bounds involving moment generating
  functions.\qed

We can now extend the Hoeffding-type lower bound (\ref{eq:32}) to
establish the asymptotic optimality of the adaptive test in
Section~\ref{sec:Mtild} that allows mid-course modification of the
maximum sample size. This adaptive test can be regarded as a
mid-course amendment of an adaptive test of $H_0: u(\theta)\le u_0$
versus $H_1: u(\theta)\ge u_1$, with a maximum sample size of $M$, to
that of $H_0$ versus $H_2: u(\theta)\ge u_2$, with a maximum sample
size of $\widetilde{M}$.  Whereas (\ref{eq:32}) provides an asymptotic
lower bound for tests of $H_0$ versus $H_1$, any test of $H_0$ versus
$H_2$ with error probabilities not exceeding $\alpha$ and
$\widetilde{\alpha}$ and taking at least $m$ and at most
$\widetilde{M}$ observations likewise satisfies
\begin{equation}\label{eq:33}E_\theta (T)\ge m\vee\left\{\widetilde{M}\wedge
    \frac{(1+o(1)) |\log\alpha|}{I_0(\theta)\vee
      I_2(\theta)}\right\}\end{equation} as $\alpha+\widetilde{\alpha}\To 0$ such that $\log\alpha\sim
\log\widetilde{\alpha}$. Note that $\Theta_1=\{\theta:u(\theta)\ge
u_1\}\subset \Theta_2=\{\theta:u(\theta)\ge u_2\}$ and therefore
$I_2(\theta)\le I_1(\theta)$. The 4-stage test in
Section~\ref{sec:Mtild}, with $M'=\widetilde{M}$, attempts to attain the asymptotic lower bound
in (\ref{eq:32}) prior to the third stage and the asymptotic lower
bound in (\ref{eq:33}) afterwards. It replaces $I_1(\theta)$ in
(\ref{eq:32}), which corresponds to early stopping for futility, by
$I_2(\theta)$ that corresponds to rejection of $H_2$ (instead of
$H_1$) in favor of $H_0$. Thus, the second-stage sample size $n_2$
corresponds to the lower bound in (\ref{eq:32}) with $\theta$
replaced by $\widehat{\theta}_m$ and $I_1$ replaced by $I_2$, while
the third-stage sample size corresponds to that in (\ref{eq:33}) with $\theta$ replaced by
$\widehat{\theta}_{n_2}$. The arguments used to prove the asymptotic
optimality of the 3-stage test in Theorem~4.1 can be readily modified
to prove the following.

\noindent\textbf{Theorem~4.2.}  \emph{Let $N^*$denote the  sample size
  of the four-stage GLR test in Section~\ref{sec:Mtild}, with
  $M'=\widetilde{M}$. Assume that $\log\alpha\sim
\log\widetilde{\alpha}$ as $\alpha+\widetilde{\alpha}\To 0$, that
(\ref{eq:17}) holds and $\widetilde{M}/|\log\alpha|\To \widetilde{A}$
with $0<a<A<\widetilde{A}$. Then $$E_\theta (N^*)\sim \left\{
  \begin{array}{ll}
m\vee (1+o(1))|\log \alpha|/I_0(\theta)&\mbox{if
  $I_0(\theta)>A^{-1}$}\\
m\vee [\widetilde{M}\wedge (1+o(1))|\log \alpha|/\{I_0(\theta)\vee
I_2(\theta)\}] &\mbox{if $I_0(\theta)<A^{-1}$.}
  \end{array}\right.$$}

The simulation results in Table~5 show that the two-stage tests of
\citet{Proschan95} and~\citet{Li02}, which use conditional power to
determine the second-stage sample size, have actual power much lower
than the adaptive test of Section~\ref{sec:three}. Without assuming a
prespecified alternative $\theta_1$, the usual approach in the
literature on sample size re-estimation considers the case $d=1$ and
$u(\theta)=\theta$ and determines the
second-stage sample size via the conditional power criterion
\begin{equation}\label{eq:19}
\widetilde n(\theta) = \min \{n \geq m : P_\theta (S_n \geq
c_{\alpha, n}|S_m) \geq 1 - \widetilde \alpha\}, \end{equation}

\noindent choosing $n_2=\widetilde n(\widehat{\theta}_m)$ if $\widehat
\theta_m \geq \theta_0 + \delta$ and stopping at the first stage due
to futility otherwise, where $P_{\theta_0}\{ S_n\ge c_{\alpha,n}\}=\alpha$ and
$\delta$ is chosen ``to set an upper bound to limit the sample size of
the second stage''; see~\citet[p.~283]{Li02}. Although the conditional
power given $\widehat \theta_m \geq \theta_0 + \delta$ is at least $1
- \widetilde \alpha$ by choosing the second-stage sample size to be
$\widetilde n (\widehat \theta_m)$, the actual (unconditional) Type II
error probability of the test at $\theta (> \theta_0)$ may
substantially exceed $\widetilde \alpha$ if $m$ is not large enough
since $P_\theta\{\widehat \theta_m < \theta_0+\delta\}$ may well
exceed $\widetilde \alpha$. Stopping due to futility at the end of the
first stage when $\widehat \theta_m < \theta_0 + \delta$ can lead to
serious loss of power of the two-stage test.  By allowing the test to
have a possible third stage, we do not have to stop prematurely when
$\widehat \theta_m$ falls below $\theta_0$, for which the conditional
power criterion (\ref{eq:19}) is not applicable. Thus, a three-stage
test that uses conditional power to determine the second-stage sample
size chooses $n_1 = m$, $n_3 = M$ and $n_2 = \min \{\widetilde n
(\widehat \theta_m), M\}$, where $\widetilde n(\theta)$ is defined by
(\ref{eq:19}) if $\theta > \theta_0$, and by
\begin{equation}\label{eq:20}\widetilde n(\theta) = \max \{m, \lceil | \log \widetilde \alpha |
  / I(\theta, \theta_1)\rceil\} \quad \mbox{if} \quad \theta \leq
  \theta_0 . \end{equation} The rejection and futility boundaries are given by (\ref{eq:6})-(\ref{eq:8}) as in
the three-stage test of Section~\ref{sec:three}. The asymptotic properties of the test, whose
sample size is 
denoted by $\widetilde{N}$, are given by the
following.

\noindent\textbf{Theorem~4.3.} \emph{Define $\eta_\theta$ for $\theta >
  \theta_0$ by
\begin{equation}\label{eq:21}\theta_0 < \eta_\theta < \theta \quad \mbox{and} \quad
  I(\eta_\theta, \theta_0) = I(\eta_\theta, \theta) . \end{equation} Then as $\alpha + \widetilde \alpha \rightarrow 0, E_\theta
(\widetilde N) \sim \widetilde n(\theta)$ if $\theta < \theta_0$ and
\begin{equation}\label{eq:22}E_\theta (\widetilde N) \sim m \vee \{M \wedge | \log \alpha | /
I(\eta_\theta, \theta_0)\} \quad \mbox{if} \quad \theta > \theta_0 .
\end{equation}}

\noindent\textit{Proof.}  Suppose $\theta>\theta_0$. From (\ref{eq:17}) and
the law of
large numbers, it follows that
$P_{\theta}\{\widehat{\theta}_m>\theta_0\}\To 1$ as $\alpha+\altild \rightarrow
0$, and therefore
\begin{equation}\label{eq:23}P_{\theta}\{\widehat{\theta}_m>\theta_0\;\mbox{and}\;
  \tilde{N}=\tilde{n}(\widehat{\theta}_m) \wedge M\} \To 1. \end{equation}
Since $\tilde{n}(\theta)$ is given by (\ref{eq:19}) in this case, application
of Theorem 2(i) of~\citet{Lai04} then yields
\begin{equation}\label{eq:24}\tilde{n}(\theta)\sim m\vee (|\log\alpha|/I(\eta_{\theta},\theta)).
\end{equation}  In view of (\ref{eq:17}), we can
apply Lebesgue's dominated convergence theorem to obtain (\ref{eq:22}) from
(\ref{eq:23}) and (\ref{eq:24}) in this case. Note in this connection that
$I(\eta_{\theta},\theta)$ is continuous in $\theta$ and that
$\widehat{\theta}_m$ converges to $\theta$ in probability.

We next consider the case $\theta<\theta_0$. Then $\tilde{n}(\theta)$
is given by (\ref{eq:20}), which is less than $M\sim
|\log\alpha|/I(\theta^*,\theta_1)$ as $\alpha+\altild \rightarrow
0$ such that $\log\alpha\sim \log\altild$. By the law of large numbers, $P_{\theta}\{\widehat{\theta}_m<\theta_0\;\mbox{and}\;
\tilde{N}=\tilde{n}(\widehat{\theta}_m)\} \To 1$. Continuity of
$I(\theta,\theta_1)$ and dominated convergence can then be used to
show that $E_{\theta}(\tilde{N})\sim \tilde{n}(\theta)$.\qed

Since $I(\eta_\theta, \theta_0) < I(\theta, \theta_0)$ by
(\ref{eq:21}), it follows from (\ref{eq:18}) and (\ref{eq:22}) that
the three-stage test using the conditional power criterion
(\ref{eq:19}) is not asymptotically efficient.
This is not surprising since (\ref{eq:19}) is the sample size for the
level-$\alpha$ FSS test to have at least $1 - \widetilde \alpha$ power
at the alternative $\theta(> \theta_0)$. However, the optimal test
with error probabilities $\alpha$ at $\theta_0$ and $\widetilde
\alpha$ at $\theta$ is Wald's sequential probability ratio test whose
expected sample size is of the smaller order $|\log \alpha | /
I(\theta, \theta_0)$ under the assumptions of Theorem~4.1.

\section{\large DISCUSSION}\label{sec:discussion}

A major drawback of the commonly used conditional power approach to
two-stage designs is pointed out in Section~\ref{sec:var_unknown}. The
actual power can be much lower than the conditional power since the
estimated alternative at the end of the first stage can be quite
different from the actual alternative. In particular, if the estimated
alternative falls in the region of the null hypothesis and misleads
one to stop for futility, there can be substantial loss of power. On
the other hand, early stopping for futility is critical for keeping
the sample size of a conditional power test within a manageable bound
$M$. Our three-stage test makes use of $M$ to come up with an implied
alternative which is used to choose the rejection and futility
boundaries appropriately so that the test does not lose much power in
comparison with the (most powerful) fixed sample size test of the null
hypothesis versus the implied alternative. This idea underlying
(\ref{eq:9})-(\ref{eq:11}) that define the stopping boundaries of
three-stage tests has been used earlier by~\citet{Lai04} to develop
efficient group sequential tests.

Our approach estimates the second-stage sample size by using an
approximation to Hoeffding's lower bound for the expected sample size
of sequential tests satisfying a prescribed Type~I error constraint and a Type~II error constraint at
the alternative that is estimated at the end of the first stage. As shown in
the simulation studies in~\citet{Bartroff08c} and in the asymptotic
theory in Section~\ref{sec:asymptotics}, this approach yields adaptive
tests that are comparable to the benchmark optimal adaptive test of
\citet{Jennison06a,Jennison06b} for a normal mean, which assumes known variance and a specified alternative. Our approach
can serve to bridge the gap between the two ``camps'' in the adaptive
design literature: One camp focuses on efficient
designs, under restrictive assumptions, that involve sufficient
statistics and optimal stopping rules, while the other camp
emphasizes flexibility to address the difficulty of coming up with realistic
alternatives at the design stage. As pointed out in
Section~\ref{sec:intro}, our approach that is built on the foundations
of sequential testing theory is able to
resolve the dilemma between efficiency and flexibility. Like the
``efficiency camp,'' it adheres to the GLR test statistics whose efficiency
is well established in the theory of FSS tests. An important
innovation is that it uses the
Markov property to compute error probabilities when the fixed sample
size is replaced by a data-dependent sample size that is based on the
estimated alternative at the end of the first stage, like the ``flexibility camp.''

\section*{\large ACKNOWLEDGMENTS}

Bartroff's work was supported by grant DMS-0403105 from the National Science
Foundation. Lai's work was supported by grants RO1-CA088890 and
DMS-0305749 from the National Institutes of Health and the National
Science Foundation.


\end{document}